\documentclass{amsart}
\usepackage{amssymb}
\usepackage{amsmath}
\usepackage{amsfonts}

\newtheorem{theorem}{Theorem}
\theoremstyle{plain}

\newtheorem{corollary}{Corollary}

\numberwithin{equation}{section}
\input{tcilatex}

\begin{document}
\title[Identities involving the $\left(h,q\right)$-Genocchi polynomials]{Identities involving 
the $\left(h,q\right)$-Genocchi polynomials and ($h,q$)-Zeta-type function}
\author{A. Bagdasaryan}
\address{Russian Academy of Sciences, Institute for Control
Sciences, 65 Profsoyuznaya, 117997 Moscow, RUSSIA}
\email{abagdasari@hotmail.com}
\author{E. \c{S}en}
\address{Department of Mathematics, Faculty of Science and Letters, 
Nam\i k Kemal University, 59030 Tekirda\u{g}, TURKEY}
\email{erdogan.math@gmail.com}
\author{Y. He}
\address{Department of Mathematics, Kunming University of Science
and Technology, Kunming, Yuannan 650500, People's Republic of China}
\email{hyyhe@yahoo.com.cn}
\author{S. Araci}
\address{Atat\"{u}rk Street, 31290 Hatay, TURKEY}
\email{mtsrkn@hotmail.com}
\author{M. Acikgoz}
\address{University of Gaziantep, Faculty of Science and Arts,
Department of Mathematics, 27310 Gaziantep, TURKEY}
\email{acikgoz@gantep.edu.tr}

\begin{abstract}
The fundamental objective of this paper is to obtain some interesting
properties for ($h,q$)-Genocchi numbers and polynomials by using the
fermionic $p$-adic $q$-integral on $\mathbb{Z}_{p}$ and mentioned in the paper $q$-Bernstein polynomials. 
By considering the $q$-Euler zeta function defined by T.~Kim, 
which can also be obtained by applying the Mellin transformation to the generating function of ($h,q$)-Genocchi
polynomials, we study ($h,q$)-Zeta-type function.
We derive symmetric properties of ($h,q$)-Zeta function and from these properties we
give symmetric property of ($h,q$)-Genocchi polynomials.

\vspace{2mm}\noindent \textsc{2010 Mathematics Subject Classification.}
11S80, 11B68.

\vspace{2mm}

\end{abstract}

\date{December 1, 2013}

\maketitle

\section{Preliminaries}

The rapid development of $q$-calculus has led to the discoveries of new
generalizations of the Bernstein polynomials and Genocchi polynomials
involving $q$-integers. The $q$-calculus theory is a novel theory that is
based on finite difference re-scaling. Remarkably, the $q$-calculus
encompasses many results of eighteenth and nineteenth century mathematics:
Euler's identities for $q$-exponential functions, Gauss's $q$-binomial
formulae, and Heine's formulae for $q$-hypergeometirc functions. Kurt Hensel
also invented $p$-adic numbers. In spite of their being already one hundred
years old, the $p$-adic numbers are still today enveloped in an aura of
mystery within the scientific community (see [1-44]).

In this paper, we also derive some interesting identities of{} $\left(
h,q\right) $-Genocchi polynomial by using $p$-adic $q$-integral on $%
\mathbb{Z}
_{p}$ and Kim's $q$-Bernstein polynomials. So, first, we list the definition
of the following notations that we use in this paper.

Let $p$ be a chosen odd prime number. Throughout this paper, the symbols $%
\mathbb{Z}
_{p},%
\mathbb{Q}
_{p},%
\mathbb{C}
,$ and $%
\mathbb{C}
_{p}$ stand for the ring of $p$-adic rational integers, the field of $p$%
-adic rational numbers, the complex number field, and the completion of
algebraic closure of $%
\mathbb{Q}
_{p}$, respectively. Let $%
\mathbb{N}
$ be the set of natural numbers and $%
\mathbb{N}
^{\ast }=%
\mathbb{N}
\cup \left\{ 0\right\} $. The $p$-adic absolute value in $%
\mathbb{C}
_{p}$ is defined by 
\begin{equation*}
\left\vert x\right\vert _{p}=p^{-r},
\end{equation*}%
where $x=p^{r}\frac{m}{n}(r\in 
\mathbb{Q}
,$ and $m,n\in 
\mathbb{Z}
$ with $\left( p,m\right) =\left( m,n\right) =\left( p,n\right) =1)$. When
one speaks of $q$-extension, $q$ is variously considered as an
indeterminate, either a complex number $q\in 
\mathbb{C}
$ or a $p$-adic number $q\in 
\mathbb{C}
_{p}$. If $q\in 
\mathbb{C}
$, we assume that $\left\vert q\right\vert <1$. If $q\in 
\mathbb{C}
_{p}$, we assume $\left\vert 1-q\right\vert _{p}<p^{-\frac{1}{p-1}}$ so that 
$q^{x}=\exp \left( x\log q\right) $ for each $\left\vert x\right\vert
_{p}\leq 1$.

The following distribution on $%
\mathbb{Z}
_{p}$ is defined by Kim as: 

for $q\in 
\mathbb{C}
_{p}$ with $\left\vert 1-q\right\vert _{p}<1$, 
\begin{equation*}
\mu _{-q}\left( x+p^{n}%
\mathbb{Z}
_{p}\right) =\left( 1+q\right) \frac{\left( -q\right) ^{x}}{\left(
1+q^{p^{n}}\right) }\text{, (for details, see \cite{Kim 5}, \cite{Kim 6}, 
\cite{Kim 7}).}
\end{equation*}

We say that $f$ is uniformly differentiable function at a point $a\in \mathbb{Z}_{p}$ and denote 
this property by the symbol $f\in UD\left(\mathbb{Z}_{p},\mathbb{C}_{p}\right) $, if the difference quotients 
\begin{equation*}
F_{f}(x,y)=\frac{f\left( x\right) -f\left( y\right) }{x-y}
\end{equation*}%
have a limit $l=f{\acute{}}\left( a\right) $ as $\left( x,y\right) \rightarrow \left( a,a\right) $.
Thus, for $f\in UD\left(\mathbb{Z}_{p},\mathbb{C}_{p}\right) $, the fermionic $p$-adic $q$-integral 
on $\mathbb{Z}_{p}$ is defined by Kim in \cite{Kim 5}, \cite{Kim 6}, \cite{Kim 7} as
follows:%
\begin{equation}
I_{-q}\left( f\right) =\int_{\mathbb{Z}_{p}}f\left( x\right) d\mu _{-q}\left( x\right) =\lim_{N\rightarrow \infty }%
\frac{1}{\left[ p^{N}\right] _{-q}}\sum_{x=0}^{p^{N}-1}\left( -1\right)
^{x}f\left( x\right) q^{x}\text{.}  \label{equation 15}
\end{equation}

So that,%
\begin{equation}
\lim_{q\rightarrow 1}I_{-q}\left( f\right) =I_{-1}\left( f\right) \text{,}
\label{equation 16}
\end{equation}%
where the notation of $I_{-1}\left( f\right) $ is called the fermionic $p$%
-adic integral on $%
\mathbb{Z}
_{p}$ (see \cite{Araci}, \cite{Araci 2}, \cite{Araci 3}, \cite{Araci 4}, 
\cite{Araci 5}, \cite{Araci 6}, \cite{Araci 8}, \cite{Kim 16}, \cite{Cangul}%
, \cite{Kim}, \cite{Kim 4}, \cite{Kim 5}, \cite{Kim 6}, \cite{Kim 7}, \cite%
{Kim 8}, \cite{Kim 10}, \cite{Kim 12}, \cite{Kim 14}, \cite{Kim 17}, \cite%
{Kim 19}, \cite{Ryoo}).

In \cite{Kim 8}, for $k,n\in 
\mathbb{N}
^{\ast }$ and $x\in \left[ 0,1\right] $, Kim's $q$-Bernstein polynomials are
defined by%
\begin{equation}
B_{k,n}\left( x,q\right) =\binom{n}{k}\left[ x\right] _{q}^{k}\left[ 1-x%
\right] _{q^{-1}}^{n-k}\text{.}  \label{equation 2}
\end{equation}

It is obvious that $\lim_{q\rightarrow 1}B_{k,n}\left( x,q\right)
=B_{k,n}\left( x\right) $ which are called the classical Bernstein
polynomials \textit{cf}. \cite{Araci 7}, \cite{Araci 1}, \cite{Araci 2}, 
\cite{Araci 4}, \cite{Araci 8}, \cite{Kim 2}, \cite{Kim 4}, \cite{Kim 8}.

In the theory of $q$-calculus for a real parameter $q\in \left( 0,1\right) $%
, $q$-analogue of $x$ is given by 
\begin{eqnarray*}
\left[ x\right] _{q} &=&\frac{1-q^{x}}{1-q}=1+q+q^{2}+...+q^{x-1}\text{,} \\
\left[ x\right] _{-q} &=&\frac{1-\left( -q\right) ^{x}}{1+q}=1+\left(
-q\right) +q^{2}+...+\left( -1\right) ^{x-1}q^{x-1}\text{.}
\end{eqnarray*}

We want to note that $\lim_{q\rightarrow 1}\left[ x\right] _{q}=x$ (see
[1-44]). Let us now take $f\left( x\right) =e^{tx}$ in (\ref{equation 16}),
then we get%
\begin{equation}
t\int_{%
\mathbb{Z}
_{p}}e^{xt}d\mu _{-1}\left( x\right) =\frac{2t}{e^{t}+1}=\sum_{n=0}^{\infty
}G_{n}\frac{t^{n}}{n!}  \label{equation 20}
\end{equation}%
where $G_{n}$ are Genocchi numbers. By using (\ref{equation 20}), we have%
\begin{equation*}
\int_{%
\mathbb{Z}
_{p}}e^{tx}d\mu _{-1}\left( x\right) =\sum_{n=0}^{\infty }\left( \frac{%
G_{n+1}}{n+1}\right) \frac{t^{n}}{n!}\text{.}
\end{equation*}

From the above, we readily see that 
\begin{equation*}
\sum_{n=0}^{\infty }\left( \int_{%
\mathbb{Z}
_{p}}x^{n}d\mu _{-1}\left( x\right) \right) \frac{t^{n}}{n!}%
=\sum_{n=0}^{\infty }\left( \frac{G_{n+1}}{n+1}\right) \frac{t^{n}}{n!}\text{%
.}
\end{equation*}

By comparing the coefficients of $\frac{t^{n}}{n!}$ in the both sides of the
above equation, we procure the following:%
\begin{equation*}
\frac{G_{n+1}}{n+1}=\int_{%
\mathbb{Z}
_{p}}x^{n}d\mu _{-1}\left( x\right) \text{, \textit{cf.} \cite{Araci 2}, 
\cite{Araci 3}, \cite{Araci 4}, \cite{Cangul}. }
\end{equation*}

In \cite{Kim 1}, the $q$-extension of Genocchi numbers are defined by%
\begin{equation*}
G_{0,q}=0,\text{ \ }q\left( qG_{q}+1\right) ^{n}+G_{n,q}=\left\{ 
\begin{array}{cc}
\left[ 2\right] _{q}, & \text{if }n=1 \\ 
0, & \text{if }n\neq 1,%
\end{array}%
\right.
\end{equation*}%
with the usual convention about replacing $\left( G_{q}\right) ^{n}$ by $%
G_{n,q}$.

Recently, several mathematicians have studied on the concept of $(h,q)$%
-Genocchi polynomials and given some new properties about these polynomials
cf. \cite{Araci 3}, \cite{Jang}, \cite{Jung}, \cite{Ryoo}. By the same
motivation, for $n\in \mathbb{N}^{\ast }$, we consider the following $\left( h,q\right) $-Genocchi
polynomials%
\begin{equation}
\frac{G_{n+1,q}^{\left( h\right) }\left( x\right) }{n+1}=\int_{%
\mathbb{Z}
_{p}}q^{\left( h-1\right) y}\left[ x+y\right] _{q}^{n}d\mu _{-q}\left(
y\right) \text{.}  \label{equation 3}
\end{equation}

In the special case $x=0$ in (\ref{equation 3}), we have $G_{n,q}^{\left(
h\right) }\left( 0\right) :=G_{n,q}^{\left( h\right) }$ that are called the $%
\left( h,q\right) $-Genocchi numbers. By (\ref{equation 3}), we derive new
relations by using aforementioned $q$-Bernstein polynomials and define their
generating function. By applying Mellin transformation to this generating
function, we obtain ($h,q$)-analogue of zeta function which interpolates ($%
h,q$)-Genocchi polynomials at negative integers. Next, we investigate
symmetric properties of the ($h,q$)-zeta function. Further, from this
investigation, we get symmetric property of ($h,q$)-Genocchi polynomials
which we present in the next sections.

\section{On the properties of the $\left( h,q\right)$-Genocchi polynomials}

By (\ref{equation 3}), we easily get%
\begin{eqnarray}
G_{n+1,q}^{\left( h\right) }\left( x\right)  &=&\left( n+1\right) \int_{%
\mathbb{Z}
_{p}}q^{\left( h-1\right) y}\left[ x+y\right] _{q}^{n}d\mu _{-q}\left(
y\right)   \notag \\
&=&\left( n+1\right) \frac{\left[ 2\right] _{q}}{\left( 1-q\right) ^{n}}%
\sum_{k=0}^{n}\binom{n}{k}q^{kx}\frac{\left( -1\right) ^{k}}{1+q^{h+k}} 
\notag \\
&=&\left[ 2\right] _{q}\left( n+1\right) \sum_{l=0}^{\infty }\left(
-1\right) ^{l}q^{hl}\left[ x+l\right] _{q}^{n}\text{.}  \label{equation 54}
\end{eqnarray}

Further, 
\begin{eqnarray*}
\sum_{n=0}^{\infty }G_{n,q}^{\left( h\right) }\left( x\right) \frac{t^{n}}{n!%
} &=&\left[ 2\right] _{q}t\sum_{n=0}^{\infty }\left( -1\right) ^{n}q^{hn}e^{%
\left[ x+n\right] _{q}t} \\
&=&\left[ 2\right] _{q}t\sum_{n=0}^{\infty }\left( -1\right)
^{n}q^{hn}e^{\left( \left[ x\right] _{q}+q^{x}\left[ n\right] _{q}\right) t}
\\
&=&\left( \frac{e^{\left[ x\right] _{q}t}}{q^{x}}\right) \left( \left[ 2%
\right] _{q}q^{x}t\sum_{n=0}^{\infty }\left( -1\right) ^{n}q^{hn}e^{\left(
q^{x}t\right) \left[ n\right] _{q}}\right) \\
&=&\left( \sum_{n=0}^{\infty }\left[ x\right] _{q}^{n}\frac{t^{n}}{n!}%
\right) \left( \sum_{n=0}^{\infty }q^{\left( n-1\right) x}G_{n,q}^{\left(
h\right) }\frac{t^{n}}{n!}\right) \\
&=&\sum_{n=0}^{\infty }\left( \sum_{k=0}^{n}\binom{n}{k}q^{\left( k-1\right)
x}\left[ x\right] _{q}^{n-k}G_{k,q}^{\left( h\right) }\right) \frac{t^{n}}{n!%
}\text{.}
\end{eqnarray*}

Therefore, we obtain the following theorem.

\begin{theorem}
For $n\in 
\mathbb{N}
^{\ast }$, we have%
\begin{equation*}
\frac{G_{n+1,q}^{\left( h\right) }\left( x\right) }{n+1}=\left[ 2\right]
_{q}\sum_{l=0}^{\infty }\left( -1\right) ^{l}q^{hl}\left[ x+l\right] _{q}^{n}%
\text{.}
\end{equation*}%
Moreover,%
\begin{equation*}
G_{n,q}^{\left( h\right) }\left( x\right) =\sum_{k=0}^{n}\binom{n}{k}%
q^{\left( k-1\right) x}G_{k,q}^{\left( h\right) }\left[ x\right]
_{q}^{n-k}=q^{-x}\left( q^{x}G_{q}^{\left( h\right) }+\left[ x\right]
_{q}\right) ^{n}
\end{equation*}%
with the usual \ convention about replacing $\left( G_{q}^{\left( h\right)
}\right) ^{n}$ by $G_{n,q}^{\left( h\right) }$.
\end{theorem}

By Theorem 1, we attain the following:%
\begin{equation}
\sum_{n=0}^{\infty }G_{n,q}^{\left( h\right) }\left( x\right) \frac{t^{n}}{n!%
}=\left[ 2\right] _{q}t\sum_{n=0}^{\infty }\left( -1\right) ^{n}q^{hn}e^{%
\left[ x+n\right] _{q}t}\text{.}  \label{equation 4}
\end{equation}

Using (\ref{equation 4}), we are now ready to obtain symmetric property of ($%
h,q $)-Genocchi polynomials, as follows: 
\begin{eqnarray*}
G_{n+1,q^{-1}}^{\left( h\right) }\left( 1-x\right) &=&(n+1)\int_{%
\mathbb{Z}
_{p}}q^{-\left( h-1\right) y}\left[ 1-x+y\right] _{q^{-1}}^{n}d\mu
_{-q^{-1}}\left( y\right) \\
&=&\left( n+1\right) \frac{\left[ 2\right] _{q^{-1}}}{\left( 1-q^{-1}\right)
^{n}}\sum_{k=0}^{n}\binom{n}{k}q^{-k\left( 1-x\right) }\left( -1\right) ^{k}%
\frac{1}{1+q^{-\left( h+k\right) }} \\
&=&\left( -1\right) ^{n}q^{h+n-1}\left( n+1\right) \frac{\left[ 2\right] _{q}%
}{\left( 1-q\right) ^{n}}\sum_{k=0}^{n}\binom{n}{k}q^{kx}\left( -1\right)
^{k}\frac{1}{1+q^{h+k}} \\
&=&\left( -1\right) ^{n}q^{h+n-1}G_{n+1,q}^{\left( h\right) }\left( x\right) 
\text{.}
\end{eqnarray*}

So, we arrive at the following theorem.

\begin{theorem}
(symmetric property of $G_{n,q}^{\left( h\right) }\left( x\right) $) Let $%
n\in 
\mathbb{N}
^{\ast }$, then we have%
\begin{equation*}
G_{n+1,q^{-1}}^{\left( h\right) }\left( 1-x\right) =\left( -1\right)
^{n}q^{h+n-1}G_{n+1,q}^{\left( h\right) }\left( x\right) \text{.}
\end{equation*}
\end{theorem}

Because of (\ref{equation 4}), we note that%
\begin{equation}
q^{h}\sum_{n=0}^{\infty }G_{n,q}^{\left( h\right) }\left( 1\right) \frac{%
t^{n}}{n!}+\sum_{n=0}^{\infty }G_{n,q}^{\left( h\right) }\frac{t^{n}}{n!}=%
\left[ 2\right] _{q}t\text{.}  \label{equation 5}
\end{equation}

By expression of (\ref{equation 5}), we derive the following recurrence
formula:%
\begin{equation}
G_{0,q}^{\left( h\right) }=0,\text{ }q^{h}G_{n,q}^{\left( h\right) }\left(
1\right) +G_{n,q}^{\left( h\right) }=\left\{ 
\begin{array}{cc}
\left[ 2\right] _{q}\text{,} & \text{if }n=1 \\ 
0\text{,} & \text{if\ }n\neq 1\text{.}%
\end{array}%
\right.  \label{equation 6}
\end{equation}

Thanks to (\ref{equation 6}) and Theorem 1, we have the following theorem.

\begin{theorem}
Let $n\in 
\mathbb{N}
^{\ast }$, then we have 
\begin{equation*}
G_{0,q}^{\left( h\right) }=0,\text{ and }q^{h-1}\left( qG_{q}^{\left(
h\right) }+1\right) ^{n}+G_{n,q}^{\left( h\right) }=\left\{ 
\begin{array}{cc}
\left[ 2\right] _{q}\text{,} & \textup{if }n=1 \\ 
0\text{,} & \textup{if }n\neq 1\text{.}%
\end{array}%
\right.
\end{equation*}%
with the usual convention about replacing $\left( G_{q}^{\left( h\right)
}\right) ^{n}$ by $G_{n,q}^{\left( h\right) }$.
\end{theorem}

For $n\in 
\mathbb{N}
$, by Theorem 1 and Theorem 3, we can proceed as follows: 
\begin{eqnarray*}
q^{2}G_{n,q}^{\left( h\right) }\left( 2\right) &=&\left( q\left(
qG_{q}^{\left( h\right) }+1\right) +1\right) ^{n} \\
&=&\sum_{k=0}^{n}\binom{n}{k}q^{k}\left( qG_{q}^{\left( h\right) }+1\right)
^{k} \\
&=&nq\left( qG_{q}^{\left( h\right) }+1\right) ^{1}+q^{1-h}\sum_{k=2}^{n}%
\binom{n}{k}q^{k}q^{h-1}\left( qG_{q}^{\left( h\right) }+1\right) ^{k} \\
&=&nq^{2-h}\left( \left[ 2\right] _{q}-\frac{\left[ 2\right] _{q}}{1+q^{h}}%
\right) -q^{1-h}\sum_{k=2}^{n}\binom{n}{k}q^{k}G_{k,q}^{\left( h\right) } \\
&=&nq^{2-h}\left[ 2\right] _{q}+q^{2-2h}G_{n,q}^{\left( h\right) }\text{ \ ,
if \ }n>1\text{.}
\end{eqnarray*}

Thus, we discover the following theorem.

\begin{theorem}
\label{thm4}Let $n\in 
\mathbb{N}
$, then we have 
\begin{equation*}
G_{n+1,q}^{\left( h\right) }\left( 2\right) =\left( n+1\right) q^{-h}\left[ 2%
\right] _{q}+q^{-2h}G_{n+1,q}^{\left( h\right) }\text{.}
\end{equation*}
\end{theorem}

With the help of Theorem 2, it is not difficult to see the following:%
\begin{eqnarray}
&&\left( n+1\right) q^{h-1}\int_{%
\mathbb{Z}
_{p}}q^{\left( h-1\right) x}\left[ 1-x\right] _{q^{-1}}^{n}d\mu _{-q}\left(
x\right)  \label{equation 7} \\
&=&q^{n+h-1}\left( -1\right) ^{n}\left( n+1\right) \int_{%
\mathbb{Z}
_{p}}q^{\left( h-1\right) x}\left[ x-1\right] _{q}^{n}d\mu _{-q}\left(
x\right)  \notag \\
&=&q^{n+h-1}\left( -1\right) ^{n}G_{n+1,q}^{\left( h\right) }\left(
-1\right) =G_{n+1,q^{-1}}^{\left( h\right) }\left( 2\right) \text{.}  \notag
\end{eqnarray}

Consequently, we state the following theorem.

\begin{theorem}
\label{thm5}The following equality holds true:%
\begin{equation*}
\left( n+1\right) q^{h-1}\int_{%
\mathbb{Z}
_{p}}q^{\left( h-1\right) x}\left[ 1-x\right] _{q^{-1}}^{n}d\mu _{-q}\left(
x\right) =G_{n+1,q^{-1}}\left( 2\right) \text{.}
\end{equation*}
\end{theorem}

On account of Theorem \ref{thm4} and Theorem \ref{thm5}, we derive the
following formula: 
\begin{eqnarray}
&&\left( n+1\right) q^{h-1}\int_{%
\mathbb{Z}
_{p}}q^{\left( h-1\right) x}\left[ 1-x\right] _{q^{-1}}^{n}d\mu _{-q}\left(
x\right)  \label{equation 8} \\
&=&\left( n+1\right) q^{h-1}\left[ 2\right] _{q}+q^{2h}G_{n+1,q^{-1}}^{%
\left( h\right) }\text{.}  \notag
\end{eqnarray}

From (\ref{equation 8}), we see that%
\begin{equation*}
\left( n+1\right) q^{h-1}\int_{%
\mathbb{Z}
_{p}}q^{\left( h-1\right) x}\left[ 1-x\right] _{q^{-1}}^{n}d\mu _{-q}=\left(
n+1\right) q^{h-1}\left[ 2\right] _{q}+q^{2h}G_{n+1,q^{-1}}^{\left( h\right)
}\text{.}
\end{equation*}

So, we deduce the following corollary.

\begin{corollary}
\label{corol1}The following equality holds true%
\begin{equation*}
\int_{%
\mathbb{Z}
_{p}}q^{\left( h-1\right) x}\left[ 1-x\right] _{q^{-1}}^{n}d\mu _{-q}\left(
x\right) =\left[ 2\right] _{q}+q^{h+1}\frac{G_{n+1,q^{-1}}^{\left( h\right) }%
}{n+1}\text{.}
\end{equation*}%

\end{corollary}

\section{New properties on the $\left( h,q\right)$-Genocchi numbers 
arising from the fermionic $p$-adic $q$-integral on $\mathbb{Z}_{p}$ 
and $q$-Bernstein polynomials}

In this part, we give some interesting relations between the $\left(
h,q\right) $-Genocchi numbers and $q$-Bernstein polynomials arising from
fermionic $p$-adic $q$-integral on $%
\mathbb{Z}
_{p}$. For $x\in 
\mathbb{Z}
_{p}$, we recall the definition of the aforementioned $q$-Bernstein
polynomials as follows:%
\begin{equation}
B_{k,n}\left( x,q\right) =\binom{n}{k}\left[ x\right] _{q}^{k}\left[ 1-x%
\right] _{q^{-1}}^{n-k}\text{, where }n,k\in 
\mathbb{N}
^{\ast }\text{.}  \label{equation 9}
\end{equation}

 By (\ref{equation 9}), Kim's $q$-Bernstein polynomials have the
following property:%
\begin{equation}
B_{k,n}\left( x,q\right) =B_{n-k,n}\left( 1-x,q^{-1}\right) \text{ (see \cite%
{Kim 8}).}  \label{equation 10}
\end{equation}

Thus, from Corollary \ref{corol1}, (\ref{equation 9}) and (\ref{equation 10}%
) we see that%
\begin{multline*}
\int_{\mathbb{Z}_{p}}B_{k,n}\left( x,q\right) q^{\left( h-1\right) x}d\mu _{-q}\left(
x\right) =\int_{\mathbb{Z}_{p}}B_{n-k,n}\left( 1-x,q^{-1}\right) q^{\left( h-1\right) x}d\mu_{-q}\left( x\right) \\
=\binom{n}{k}\sum_{l=0}^{k}\binom{k}{l}\left( -1\right) ^{k+l}\int_{\mathbb{Z}_{p}}q^{\left( h-1\right) x}\left[ 1-x\right] _{q^{-1}}^{n-l}d\mu
_{-q}\left( x\right) \\
=\binom{n}{k}\sum_{l=0}^{k}\binom{k}{l}\left( -1\right) ^{k+l}\left( \left[
2\right] _{q}+\frac{q^{h+1}}{n-l+1}G_{n+1,q^{-1}}^{\left( h\right) }\right) .
\end{multline*}

For $n$, $k\in 
\mathbb{N}
^{\ast }$ with $n>k$, we compute%
\begin{eqnarray}
&&\int_{%
\mathbb{Z}
_{p}}B_{k,n}\left( x,q\right) q^{\left( h-1\right) x}d\mu _{-q}\left(
x\right)  \notag \\
&=&\binom{n}{k}\sum_{l=0}^{k}\binom{k}{l}\left( -1\right) ^{k+l}\left( \left[
2\right] _{q}+\frac{q^{h+1}}{n-l+1}G_{n+1,q^{-1}}^{\left( h\right) }\right)
\label{equation 11} \\
&=&\left\{ 
\begin{array}{cc}
\left[ 2\right] _{q}+\frac{q^{h+1}}{n+1}G_{n+1,q^{-1}}^{\left( h\right) } & 
\text{if }k=0, \\  \\
\binom{n}{k}\sum_{l=0}^{k}\binom{k}{l}\left( -1\right) ^{k+l}\left( \left[ 2%
\right] _{q}+\frac{q^{h+1}}{n-l+1}G_{n+1,q^{-1}}^{\left( h\right) }\right) & 
\text{if }k>0.%
\end{array}%
\right.  \notag
\end{eqnarray}

Let us take the fermionic $p$-adic $q$-integral on $%
\mathbb{Z}
_{p}$ for the $q$-Bernstein polynomials of degree $n$ as follows:%
\begin{gather}
\int_{%
\mathbb{Z}
_{p}}B_{k,n}\left( x,q\right) q^{\left( h-1\right) x}d\mu _{-q}\left(
x\right) =\binom{n}{k}\int_{%
\mathbb{Z}
_{p}}\left[ x\right] _{q}^{k}\left[ 1-x\right] _{q^{-1}}^{n-k}q^{\left(
h-1\right) x}d\mu _{-q}\left( x\right)  \label{equation 12} \\
=\binom{n}{k}\sum_{l=0}^{n-k}\binom{n-k}{l}\left( -1\right) ^{l}\frac{%
G_{l+k+1,q}^{\left( h\right) }}{l+k+1}.  \notag
\end{gather}

Therefore, by (\ref{equation 11}) and (\ref{equation 12}), we attain the
following theorem.

\begin{theorem}
Let $n,k\in 
\mathbb{N}
^{\ast }$ with $n>k$. Then we have%
\begin{multline*}
\sum_{l=0}^{n-k}\binom{n-k}{l}\left( -1\right) ^{l}\frac{G_{l+k+1,q}^{\left(
h\right) }}{l+k+1} \\
=\left\{ 
\begin{array}{cc}
\left[ 2\right] _{q}+\frac{q^{h+1}}{n+1}G_{n+1,q^{-1}}^{\left( h\right) } & 
\textup{if }k=0, \\  \\
\sum_{l=0}^{k}\binom{k}{l}\left( -1\right) ^{k+l}\left( \left[ 2\right] _{q}+%
\frac{q^{h+1}}{n-l+1}G_{n+1,q^{-1}}^{\left( h\right) }\right) & \textup{if }%
k>0.%
\end{array}%
\right.
\end{multline*}
\end{theorem}

Putting $k=0$ in the above theorem, we procure the following corollary.

\begin{corollary}
The following holds true:%
\begin{equation*}
\sum_{l=0}^{n}\binom{n}{l}\left( -1\right) ^{l}\frac{G_{l+1,q}^{\left(
h\right) }}{l+1}=\left[ 2\right] _{q}+\frac{q^{h+1}}{n+1}G_{n+1,q^{-1}}^{%
\left( h\right) }\text{.}
\end{equation*}
\end{corollary}

Let $n_{1},n_{2},k\in 
\mathbb{N}
^{\ast }$ with $n_{1}+n_{2}>2k$. Then we get%
\begin{multline*}
\int_{\mathbb{Z}_{p}}B_{k,n_{1}}\left( x,q\right) B_{k,n_{2}}\left( x,q\right) q^{\left(
h-1\right) x}d\mu _{-q}\left( x\right)   \\
=\binom{n_{1}}{k}\binom{n_{2}}{k}\sum_{l=0}^{2k}\binom{2k}{l}\left(
-1\right) ^{2k+l}\int_{\mathbb{Z}_{p}}\left[ 1-x\right] _{q^{-1}}^{n_{1}+n_{2}-l}q^{\left( h-1\right) x}d\mu
_{-q}\left( x\right)  \notag \\
=\binom{n_{1}}{k}\binom{n_{2}}{k}\sum_{l=0}^{2k}\binom{2k}{l}\left(
-1\right) ^{2k+l}\left( \left[ 2\right] _{q}+\frac{q^{h+1}}{n_{1}+n_{2}-l+1}%
G_{n_{1}+n_{2}-l+1,q^{-1}}^{\left( h\right) }\right)   \\
=\left\{ 
\begin{array}{cc}
\left[ 2\right] _{q}+\frac{q^{h+1}}{n_{1}+n_{2}+1}G_{n_{1}+n_{2}+1,q^{-1}}^{%
\left( h\right) } & \textup{if }k=0, \\  \\
\binom{n_{1}}{k}\binom{n_{2}}{k}\sum_{l=0}^{2k}\binom{2k}{l}\left( -1\right)
^{2k+l}\left( \left[ 2\right] _{q}+\frac{q^{h+1}}{n_{1}+n_{2}-l+1}%
G_{n_{1}+n_{2}-l+1,q^{-1}}^{\left( h\right) }\right) & \textup{if }k\neq 0.%
\end{array}%
\right.  
\end{multline*}

As a result, we state the following theorem.

\begin{theorem}
Let $n_{1},n_{2},k\in 
\mathbb{N}
^{\ast }$ with $n_{1}+n_{2}>2k,$ then we get%
\begin{multline*}
\int_{\mathbb{Z}_{p}}B_{k,n_{1}}\left( x,q\right) B_{k,n_{2}}\left( x,q\right) q^{\left(
h-1\right) x}d\mu _{-q}\left( x\right) \\
=\left\{ 
\begin{array}{cc}
\left[ 2\right] _{q}+\frac{q^{h+1}}{n_{1}+n_{2}+1}G_{n_{1}+n_{2}+1,q^{-1}}^{%
\left( h\right) } & \textup{if }k=0, \\  \\
\binom{n_{1}}{k}\binom{n_{2}}{k}\sum_{l=0}^{2k}\binom{2k}{l}\left( -1\right)
^{2k+l}\left( \left[ 2\right] _{q}+\frac{q^{h+1}}{n_{1}+n_{2}-l+1}%
G_{n_{1}+n_{2}-l+1,q^{-1}}^{\left( h\right) }\right) & \textup{if }k\neq 0.%
\end{array}%
\right.
\end{multline*}
\end{theorem}

From the binomial theorem, we can derive the following equation.%
\begin{multline}
\int_{%
\mathbb{Z}
_{p}}B_{k,n_{1}}\left( x,q\right) B_{k,n_{2}}\left( x,q\right) q^{\left(
h-1\right) x}d\mu _{-q}\left( x\right) \label{equation 14}  \\
=\displaystyle\prod\limits_{i=1}^{2}\binom{n_{i}}{k}%
\sum_{l=0}^{n_{1}+n_{2}-2k}\binom{n_{1}+n_{2}-2k}{l}\left( -1\right)
^{l}\int_{%
\mathbb{Z}
_{p}}\left[ x\right] _{q}^{2k+l}q^{\left( h-1\right) x}d\mu _{-q}\left(
x\right)   \\
=\displaystyle\prod\limits_{i=1}^{2}\binom{n_{i}}{k}%
\sum_{l=0}^{n_{1}+n_{2}-2k}\binom{n_{1}+n_{2}-2k}{l}\left( -1\right) ^{l}%
\frac{G_{l+2k+1,q}^{\left( h\right) }}{l+2k+1}.  
\end{multline}

Thus, by (\ref{equation 14}), we obtain the following theorem.

\begin{theorem}
Let $n_{1},n_{2},k\in 
\mathbb{N}
^{\ast }$ with $n_{1}+n_{2}>2k,$ then we have%
\begin{multline*}
\sum_{l=0}^{n_{1}+n_{2}-2k}\binom{n_{1}+n_{2}-2k}{l}\left( -1\right) ^{l}%
\frac{G_{l+2k+1,q}^{\left( h\right) }}{l+2k+1} \\
=\left\{ 
\begin{array}{cc}
\left[ 2\right] _{q}+\frac{q^{h+1}}{n_{1}+n_{2}+1}G_{n_{1}+n_{2}+1,q^{-1}}^{%
\left( h\right) } & \textup{if }k=0, \\  \\
\sum_{l=0}^{2k}\binom{2k}{l}\left( -1\right) ^{2k+l}\left( \left[ 2\right]
_{q}+\frac{q^{h+1}}{n_{1}+n_{2}-l+1}G_{n_{1}+n_{2}-l+1,q^{-1}}^{\left(
h\right) }\right) & \textup{if }k\neq 0.%
\end{array}%
\right.
\end{multline*}
\end{theorem}

Now also, by the same method, substituting $k=0$ in the above theorem, we
discover the following corollary.

\begin{corollary}
The following identity%
\begin{equation*}
\sum_{l=0}^{n_{1}+n_{2}}\binom{n_{1}+n_{2}}{l}\left( -1\right) ^{l}\frac{%
G_{l+1,q}^{\left( h\right) }}{l+1}=\left[ 2\right] _{q}+q^{h+1}\frac{%
G_{n_{1}+n_{2}+1,q^{-1}}^{\left( h\right) }}{n_{1}+n_{2}+1}
\end{equation*}%
holds true.
\end{corollary}

For $x\in 
\mathbb{Z}
_{p}$ and $s\in 
\mathbb{N}
$ with $s\geq 2,$ let $n_{1},n_{2},...,n_{s},k\in 
\mathbb{N}
^{\ast }$ with $\sum_{l=1}^{s}n_{l}>sk$. Then we take the fermionic $p$-adic $q$-integral on 
$\mathbb{Z}_{p}$ for the $q$-Bernstein polynomials of degree $n$ as follows:%
\begin{multline*}
\int_{\mathbb{Z}_{p}}\underset{\text{s-times}}{\underbrace{B_{k,n_{1}}\left( x,q\right)
B_{k,n_{2}}\left( x,q\right) ...B_{k,n_{s}}\left( x,q\right) }q^{\left(
h-1\right) x}}d\mu _{-q}\left( x\right) \\
=\displaystyle\prod\limits_{i=1}^{s}\binom{n_{i}}{k}\int_{\mathbb{Z}_{p}}\left[ x\right] _{q}^{sk}\left[ 1-x\right]
_{q^{-1}}^{n_{1}+n_{2}+...+n_{s}-sk}q^{\left( h-1\right) x}d\mu _{-q}\left(
x\right) \\
=\displaystyle\prod\limits_{i=1}^{s}\binom{n_{i}}{k}\sum_{l=0}^{sk}\binom{%
sk}{l}\left( -1\right) ^{l+sk}\int_{\mathbb{Z}_{p}}\left[ 1-x\right] _{q^{-1}}^{n_{1}+n_{2}+...+n_{s}-sk}q^{\left(
h-1\right) x}d\mu _{-q}\left( x\right) \\
=\left\{ 
\begin{array}{cc}
\left[ 2\right] _{q}+\frac{q^{h+1}}{n_{1}+n_{2}+...+n_{s}+1}%
G_{n_{1}+n_{2}+...+n_{s}+1,q^{-1}}^{\left( h\right) } & \text{if }k=0\text{,}
\\  \\
\displaystyle\prod\limits_{i=1}^{s}\binom{n_{i}}{k}\sum_{l=0}^{sk}\binom{sk}{%
l}\left( -1\right) ^{sk+l}  \times \\
\times \left( \left[ 2\right] _{q}+\frac{q^{h+1}}{%
n_{1}+n_{2}+...+n_{s}-l+1}G_{n_{1}+n_{2}+...+n_{s}-l+1,q^{-1}}^{\left(
h\right) }\right) & \text{if }k\neq 0\text{.}%
\end{array}%
\right.
\end{multline*}

Consequently, we obtain the following theorem.

\begin{theorem}
Let $s\in 
\mathbb{N}
$ with $s\geq 2$, let $n_{1},n_{2},...,n_{s},k\in 
\mathbb{N}
^{\ast }$ with $\sum_{l=1}^{s}n_{l}>sk$. Then we have%
\begin{multline*}
\int_{\mathbb{Z}_{p}}\left( \displaystyle\prod\limits_{i=1}^{s}B_{k,n_{i}}\left( x\right)
\right) q^{\left( h-1\right) x}d\mu _{-q}\left( x\right) \\
=\left\{ 
\begin{array}{cc}
\left[ 2\right] _{q}+\frac{q^{h+1}}{n_{1}+n_{2}+...+n_{s}+1}%
G_{n_{1}+n_{2}+...+n_{s}+1,q^{-1}}^{\left( h\right) } & \textup{if }k=0%
\text{,} \\  \\
\displaystyle\prod\limits_{i=1}^{s}\binom{n_{i}}{k}\sum_{l=0}^{sk}\binom{sk}{%
l}\left( -1\right) ^{sk+l} \times \\
\times \left( \left[ 2\right] _{q}+\frac{q^{h+1}}{%
n_{1}+n_{2}+...+n_{s}-l+1}G_{n_{1}+n_{2}+...+n_{s}-l+1,q^{-1}}^{\left(
h\right) }\right) & \textup{if }k\neq 0\text{.}%
\end{array}%
\right.
\end{multline*}
\end{theorem}

From the definition of $q$-Bernstein polynomials and the binomial theorem,
we easily see that%
\begin{multline}
\int_{%
\mathbb{Z}
_{p}}\underset{s-times}{\underbrace{B_{k,n_{1}}\left( x,q\right)
B_{k,n_{2}}\left( x,q\right) ...B_{k,n_{s}}\left( x,q\right) }q^{\left(
h-1\right) x}}d\mu _{-q}\left( x\right)  \label{equation 17} \\
=\displaystyle\prod\limits_{i=1}^{s}\binom{n_{i}}{k}%
\sum_{l=0}^{n_{1}+...+n_{s}-sk}\binom{\sum_{d=1}^{s}\left( n_{d}-k\right) }{l%
}\left( -1\right) ^{l}\int_{%
\mathbb{Z}
_{p}}\left[ x\right] _{q}^{sk+l}q^{\left( h-1\right) x}d\mu _{-q}\left(
x\right)   \\
=\displaystyle\prod\limits_{i=1}^{s}\binom{n_{i}}{k}%
\sum_{l=1}^{n_{1}+...+n_{s}-sk}\binom{\sum_{d=1}^{s}\left( n_{d}-k\right) }{l%
}\left( -1\right) ^{l}\frac{G_{l+sk+1,q}^{\left( h\right) }}{l+sk+1}.  
\end{multline}

Thus, from (\ref{equation 17}), we discover the following theorem.

\begin{theorem}
Let $s\in 
\mathbb{N}
$ with $s\geq 2$, let $n_{1},n_{2},...,n_{s},k\in 
\mathbb{N}
^{\ast }$ with $\sum_{l=1}^{s}n_{l}>sk$, then the following identity holds
true: 
\begin{multline*}
\sum_{l=0}^{n_{1}+...+n_{s}-sk}\binom{\sum_{d=1}^{s}\left( n_{d}-k\right) 
}{l}\left( -1\right) ^{l}\frac{G_{l+sk+1,q}^{\left( h\right) }}{l+sk+1} \\  
=\left\{ 
\begin{array}{cc}
\left[ 2\right] _{q}+\frac{q^{h+1}}{n_{1}+n_{2}+...+n_{s}+1}%
G_{n_{1}+n_{2}+...+n_{s}+1,q^{-1}}^{\left( h\right) } & \textup{if }k=0%
\text{,} \\ \\
\sum_{l=0}^{sk}\binom{sk}{l}\left( -1\right) ^{sk+l} \times \\ 
\times  \left( \left[ 2\right]
_{q}+\frac{q^{h+1}}{n_{1}+n_{2}+...+n_{s}-l+1}%
G_{n_{1}+n_{2}+...+n_{s}-l+1,q^{-1}}^{\left( h\right) }\right) & \textup{if 
}k\neq 0\text{.}%
\end{array}%
\right.
\end{multline*}
\end{theorem}

Taking $k=0$ in the above theorem, we deduce the following.

\begin{corollary}
The identity%
\begin{equation*}
\sum_{l=0}^{n_{1}+...+n_{s}}\binom{n_{1}+...+n_{s}}{l}\left( -1\right) ^{l}%
\frac{G_{l+1,q}^{\left( h\right) }}{l+1}=\left[ 2\right] _{q}+q^{h+1}\frac{%
G_{n_{1}+n_{2}+...+n_{s}+1,q^{-1}}^{\left( h\right) }}{%
n_{1}+n_{2}+...+n_{s}+1}
\end{equation*}%
holds true.
\end{corollary}

\section{\textbf{Further Remarks}}

In this final part, we consider the $q$-Euler Zeta function in $%
\mathbb{C}
$, which is defined by Kim \cite{Kim 5}%
\begin{equation}
\zeta _{q}^{\left( h\right) }\left( s,x\right) =\left[ 2\right]
_{q}\sum_{m=0}^{\infty }\frac{\left( -1\right) ^{m}q^{mh}}{\left[ m+x\right]
_{q}^{s}}  \label{equation 49}
\end{equation}%
where $q\in 
\mathbb{C}
$, $h\in 
\mathbb{N}
$ and $\Re e\left( s\right) >1$. It is clear that for the special case $h=0$
and $q\rightarrow 1$ in (\ref{equation 49}), it reduces to the ordinary
Hurwitz-Euler zeta function. Now, we consider (\ref{equation 49}) in the form%
\begin{equation*}
\zeta _{q^{a}}^{\left( h\right) }\left( s,bx+\frac{bj}{a}\right) =\left[ 2%
\right] _{q^{a}}\sum_{m=0}^{\infty }\frac{\left( -1\right) ^{m}q^{mah}}{%
\left[ m+bx+\frac{bj}{a}\right] _{q^{a}}^{s}}\text{.}
\end{equation*}

By using some operations to the above identity, that is, for any positive
integers $m$ and $b$, there exist unique non-negative integers $k$ and $i$
such that $m=bk+i$ with $0\leq i\leq b-1$ and thus for $a\equiv 1(\func{mod}%
2)$ and $b\equiv 1(\func{mod}2)$, we can compute as follows:%
\begin{align}
\zeta _{q^{a}}^{\left( h\right) }\left( s,bx+\frac{bj}{a}\right) & =\left[ a%
\right] _{q}^{s}\left[ 2\right] _{q^{a}}\sum_{m=0}^{\infty }\frac{\left(
-1\right) ^{m}q^{mah}}{\left[ ma+abx+bj\right] _{q^{a}}^{s}}
\label{equation 50} \\
& =\left[ a\right] _{q}^{s}\left[ 2\right] _{q^{a}}\sum_{m=0}^{\infty
}\sum_{i=0}^{b-1}\frac{\left( -1\right) ^{i+mb}q^{\left( i+mb\right) ah}}{%
\left[ \left( i+mb\right) a+abx+bj\right] _{q^{a}}^{s}}  \notag \\
& =\left[ a\right] _{q}^{s}\left[ 2\right] _{q^{a}}\sum_{i=0}^{b-1}\left(
-1\right) ^{i}q^{iah}\sum_{m=0}^{\infty }\frac{\left( -1\right) ^{m}q^{mbah}%
}{\left[ ab\left( m+x\right) +ai+bj\right] _{q^{a}}^{s}}\text{.}  \notag
\end{align}

From this, we see that%
\begin{gather}
\sum_{j=0}^{a-1}\left( -1\right) ^{j}q^{jbh}\zeta _{q^{a}}^{\left( h\right)
}\left( s,bx+\frac{bj}{a}\right) =  \label{equation 51} \\
\left[ a\right] _{q}^{s}\left[ 2\right] _{q^{a}}\sum_{j=0}^{a-1}\left(
-1\right) ^{j}q^{jbh}\sum_{i=0}^{b-1}\left( -1\right)
^{i}q^{iah}\sum_{m=0}^{\infty }\frac{\left( -1\right) ^{m}q^{mbah}}{\left[
ab\left( m+x\right) +ai+bj\right] _{q}^{s}}\text{.}  \notag
\end{gather}

Replacing $a$ by $b$ and $j$ by $i$ in (\ref{equation 50}), we derive the
following%
\begin{equation*}
\zeta _{q^{b}}^{\left( h\right) }\left( s,ax+\frac{ai}{b}\right) =\left[ b%
\right] _{q}^{s}\left[ 2\right] _{q^{b}}\sum_{j=0}^{a-1}\left( -1\right)
^{j}q^{jbh}\sum_{m=0}^{\infty }\frac{\left( -1\right) ^{m}q^{mbah}}{\left[
ab\left( m+x\right) +ai+bj\right] _{q}^{s}}\text{.}
\end{equation*}

By considering the above identity in (\ref{equation 51}), we can easily
state the following theorem.

\begin{theorem}
\label{thm}The following holds true%
\begin{equation*}
\frac{\left[ 2\right] _{q^{b}}}{\left[ a\right] _{q}^{s}}\sum_{i=0}^{a-1}%
\left( -1\right) ^{i}q^{ibh}\zeta _{q^{a}}^{\left( h\right) }\left( s,bx+%
\frac{bi}{a}\right) =\frac{\left[ 2\right] _{q^{a}}}{\left[ b\right] _{q}^{s}%
}\sum_{i=0}^{b-1}\left( -1\right) ^{i}q^{iah}\zeta _{q^{b}}^{\left( h\right)
}\left( s,ax+\frac{ai}{b}\right) \text{.}
\end{equation*}
\end{theorem}

Now, setting $b=1$ in the above theorem, we easily procure the following
distribution formula%
\begin{equation}
\zeta _{q}^{\left( h\right) }\left( s,ax\right) =\frac{\left[ 2\right] _{q}}{%
\left[ 2\right] _{q^{a}}\left[ a\right] _{q}^{s}}\sum_{i=0}^{a-1}\left(
-1\right) ^{i}q^{ih}\zeta _{q^{a}}^{\left( h\right) }\left( s,x+\frac{i}{a}%
\right) \text{.}  \label{equation 52}
\end{equation}

Putting $a=2$ in (\ref{equation 52}), it leads to the following corollary.

\begin{corollary}
The following identity holds true:%
\begin{equation*}
\zeta _{q}^{\left( h\right) }\left( s,2x\right) =\frac{\left[ 2\right] _{q}}{%
\left[ 2\right] _{q^{2}}\left[ 2\right] _{q^{\alpha }}^{s}}\left( \zeta
_{q^{2}}^{\left( h\right) }\left( s,x\right) -q^{h}\zeta _{q^{2}}^{\left(
h\right) }\left( s,x+\frac{1}{2}\right) \right) \text{.}
\end{equation*}
\end{corollary}

By (\ref{equation 54}) and (\ref{equation 49}), we have%
\begin{equation}
\frac{G_{n+1,q}^{\left( h\right) }\left( x\right) }{n+1}=\zeta _{q}^{\left(
h\right) }\left( -n,x\right) \text{.}  \label{equation 55}
\end{equation}

By Theorem \ref{thm} and (\ref{equation 55}), we have the symmetric property
of ($h,q$)-Genocchi polynomials by the following theorem.

\begin{theorem}
The following identity%
\begin{multline*}
\left[ 2\right] _{q^{b}}\left[ a\right] _{q}^{m-1}\sum_{j=0}^{a-1}\left(
-1\right) ^{i}q^{ibh}G_{m,q^{a}}^{\left( h\right) }\left( bx+\frac{bi}{a}%
\right) \\
=\left[ 2\right] _{q^{a}}\left[ b\right] _{q}^{m-1}\sum_{i=0}^{b-1}%
\left( -1\right) ^{i}q^{iah}G_{m,q^{b}}^{\left( h\right) }\left( ax+\frac{ai%
}{b}\right)
\end{multline*}%
holds true.
\end{theorem}

Now also, setting $b=1$ and replacing $x$ by $\frac{x}{a}$ in the above
theorem, we can rewrite the following ($h,q$)-Genocchi polynomials.%
\begin{equation*}
G_{n,q}^{\left( h\right) }\left( x\right) =\frac{\left[ 2\right] _{q}}{\left[
2\right] _{q^{a}}}\left[ a\right] _{q}^{n-1}\sum_{i=0}^{a-1}\left( -1\right)
^{i}q^{ih}G_{n,q^{a}}^{\left( h\right) }\left( \frac{x+i}{a}\right) \text{ }%
\left( 2\nmid a\right) \text{.}
\end{equation*}

On account of (\ref{equation 4}), we develop as follows: 
\begin{eqnarray}
\sum_{n=0}^{\infty }G_{n,q}^{\left( h\right) }\left( x+y\right) \frac{t^{n}}{%
n!} &=&\left[ 2\right] _{q}t\sum_{m=0}^{\infty }\left( -1\right)
^{m}q^{mh}e^{t\left[ x+y+m\right] _{q}}  \label{equation 19} \\
&=&\left[ 2\right] _{q}t\sum_{m=0}^{\infty }\left( -1\right) ^{m}q^{mh}e^{t%
\left[ y\right] _{q}}e^{\left( q^{y}t\right) \left[ x+m\right] _{q}}  \notag
\\
&=&\left( \sum_{n=0}^{\infty }\left[ y\right] _{q}^{n}\frac{t^{n}}{n!}%
\right) \left( \sum_{n=0}^{\infty }q^{\left( n-1\right) y}G_{n,q}^{\left(
h\right) }\left( x\right) \frac{t^{n}}{n!}\right) \text{.}
\label{equation 18}
\end{eqnarray}

By using Cauchy product in (\ref{equation 18}), we see that (\ref{equation 18}) can be written as %
\begin{equation*}
\sum_{n=0}^{\infty }\left( \sum_{j=0}^{n}\binom{n}{j}q^{\left( j-1\right)
y}G_{j,q}^{\left( h\right) }\left( x\right) \left[ y\right] _{q^{\alpha
}}^{n-j}\right) \frac{t^{n}}{n!}\text{.}
\end{equation*}

Thus, by comparing the coefficients of $\frac{t^{n}}{n!}$ in (\ref{equation
19}) and (\ref{equation 18}), we state the following corollary.

\begin{corollary}
The following equality holds true:%
\begin{equation}
G_{n,q}^{\left( h\right) }\left( x+y\right) =\sum_{j=0}^{n}\binom{n}{j}%
q^{\left( j-1\right) y}G_{j,q}^{\left( h\right) }\left( x\right) \left[ y%
\right] _{q}^{n-j}\text{.}  \label{equation 53}
\end{equation}
\end{corollary}

By using (\ref{equation 53}), after some computations, we readily derive the
following symmetric relation:

\begin{theorem}
The following equality holds:%
\begin{multline*}
\left[ 2\right] _{q^{b}}\sum_{i=0}^{m}\binom{m}{i}\left[ a\right] _{q}^{i-1}%
\left[ b\right] _{q}^{m-i}G_{i,q^{a}}^{\left( h\right) }\left( bx\right)
S_{m-i;q^{b}}^{\left( h+i-1\right) }\left( a\right) \\
=\left[ 2\right] _{q^{a}}\sum_{i=0}^{m}\binom{m}{i}\left[ b\right] _{q}^{i-1}%
\left[ a\right] _{q}^{m-i}G_{i,q^{b}}^{\left( h\right) }\left( ax\right)
S_{m-i;q^{a}}^{\left( h+i-1\right) }\left( b\right)
\end{multline*}%
where $S_{m;q}^{\left( i\right) }\left( a\right) =\sum_{j=0}^{a-1}\left(
-1\right) ^{j}q^{ji}\left[ j\right] _{q}^{m}$.
\end{theorem}

\end{document}